\documentclass[12pt]{amsart}

\usepackage{amsmath,xspace,amssymb,mathrsfs}
\usepackage{color}

\input xy
\xyoption{all}
\xyoption{2cell}
\UseAllTwocells
\CompileMatrices

\newcommand{\Spec}{\operatorname{Spec}}
\renewcommand{\phi}{\varphi}

\newcommand{\Ker}{\operatorname{Ker}}
\newcommand{\Ima}{\operatorname{Im}}
\newcommand{\Min}{\operatorname{Min}}

\newcommand{\V}{\operatorname{V}}
\newcommand{\Max}{\operatorname{Max}}

\newtheorem{proposition}{Proposition}[section]
\newtheorem{lemma}[proposition]{Lemma} 
\newtheorem{corollary}[proposition]{Corollary}
\newtheorem{theorem}[proposition]{Theorem}

\newtheorem{prop-def}[proposition]{Proposition and definition}

  \theoremstyle{definition}

  \newtheorem{remark}[proposition]{Remark}

\begin{document}

\title[Zariski compactness of minimal spectrum]{Zariski compactness of minimal spectrum and flat compactness of maximal spectrum}

\author[A. Tarizadeh]{Abolfazl Tarizadeh}
\address{ Department of Mathematics, Faculty of Basic Sciences, University of Maragheh \\
P. O. Box 55136-553, Maragheh, Iran.
 }
\email{ebulfez1978@gmail.com}


\begin{abstract} In this article, Zariski compactness of the minimal spectrum and flat compactness of the maximal spectrum are characterized.

\end{abstract}

\maketitle

\section{Introduction}

The minimal spectrum of a commutative ring, specially its compactness, has been the main topic of many articles in the literature over the years and it is still of current interest, see e.g. \cite{Artico-Marconi}, \cite{Henriksen-Jerison}, \cite{Hochster}, \cite{Hong et all}, \cite{Kist}, \cite{Knox et al}, \cite{Matlis}, \cite{Mewborn}, \cite{Quentel}, \cite{Schwartz-Tressl}. Amongst them, the well-known result of  Quentel \cite[Proposition 1]{Quentel} can be considered as one of the most important results in this context. But his proof, as presented there, is  sketchy. In fact, it is merely a plan of the proof, not the proof itself. In the present article,  i.e. Theorem \ref{th 2} and Corollary \ref{coro 2}, a new and purely algebraic proof is given for this non-trivial result. Dually, a new result is also given for the compactness of the maximal spectrum with respect to the flat topology, see Theorem \ref{th 1}.\\

In Theorem \ref{lem 3}, the patch closures are computed in a certain way. This result plays a major role in proving Theorem \ref{th 2}. The noetherianess of the prime spectrum with respect to the Zariski topology is also characterized, see Theorem \ref{th 3}. It is worth mentioning that in \cite{Ohm} and \cite{Rush}, one can find other nontrivial characterizations of noetherianness of the prime spectrum. \\

\section{Preliminaries}

Here, we briefly recall some material which is needed in the sequel.\\

In this paper, all of the rings are commutative. Every ring map $\phi:A\rightarrow B$ induces a map  $\phi^{\ast}:\Spec(B)\rightarrow\Spec(A)$ between the corresponding prime spectra which maps each prime ideal $\mathfrak{p}$ of $B$ into $\phi^{-1}(\mathfrak{p})$.\\

Mel Hochster in his seminal work \cite[Proposition 8]{Hochster}, discovered a new spectral topology on the underlying set of a given spectral space which behaves as the dual of the original topology. It is called the flat (or, inverse) topology. In this paper, we need to express explicitly the flat topology on the prime spectrum. Thus, let $R$ be a commutative ring. Then the collection of subsets $V(f)=\{\mathfrak{p}\in\Spec R: f\in\mathfrak{p}\}$ with $f\in R$ forms a sub-basis for the opens of the flat topology over $\Spec R$. Moreover, there is a (unique) topology over $\Spec R$ such that the collection of subsets $D(f)\cap V(g)$ with $f,g\in R$ forms a sub-basis for the opens of this topology. It is called the patch (or, constructible) topology. The patch topology is finer than the Zariski and flat topologies. Obviously, every map $\Spec B\rightarrow\Spec A$ induced by a ring map $A\rightarrow B$ is continuous with respect to the flat topology. It is also continuous with respect to the patch topology. The flat topology behaves as the dual of the Zariski topology. For instance, if $\mathfrak{p}$ is a prime ideal of $R$, then its closure with respect to the flat topology comes from the canonical ring map $R\rightarrow R_{\mathfrak{p}}$. In fact, $\Lambda(\mathfrak{p})=\{\mathfrak{q}\in\Spec R: \mathfrak{q}\subseteq\mathfrak{p}\}$. Here, $\Lambda(\mathfrak{p})$ denotes the closure of $\{\mathfrak{p}\}$ in $\Spec R$ with respect to the flat topology. By contrast, the Zariski closure of this point comes from the canonical ring map $R\rightarrow R/\mathfrak{p}$.  \\

It is well known that the Zariski closed subsets of $\Spec(R)$ are precisely of the form $\Ima\pi^{\ast}$ where $\pi:R\rightarrow R/I$ is the canonical ring map. One can show that the patch closed subsets of $\Spec(R)$ are precisely of the form $\Ima\phi^{\ast}$ where $\phi:R\rightarrow A$ is a ring map. Moreover, the flat closed subsets of $\Spec(R)$ are precisely of the form $\Ima\phi^{\ast}$ where $\phi:R\rightarrow A$ is a flat ring map. \\

A ring $R$ is said to be absolutely flat (or, von-Neumann regular) if
each $R-$module is $R-$flat. This is equivalent to the statement that each element $f\in R$ can be written as $f=f^{2}g$ for some $g\in R$. Clearly, absolutely flat rings are stable under taking quotients and localizations. A direct product of rings $(A_{i})$ is absolutely flat if and only if each $A_{i}$ is absolutely flat. Every prime ideal of an absolutely flat ring is maximal since an absolutely flat domain is a field. Every absolutely flat ring $R$ is reduced because for each $f\in R$ and for each natural number $n\geq2$ there exists some $g\in R$ such that $f=f^{n}g^{n-1}$. A ring $R$ is absolutely flat if and only if for each maximal ideal $\mathfrak{p}$ of $R$, $R_{\mathfrak{p}}$ is absolutely flat. \\

Let $S$ be a subset of a ring $R$. The polynomial ring $R[x_{s} : s\in S]$ modulo $I$ is denoted by $S^{(-1)}R$ where the ideal $I$ is generated by elements of the form $sx_{s}^{2}-x_{s}$ and $s^{2}x_{s}-s$ with $s\in S$. The ring $S^{(-1)}R$ is called the pointwise localization of $R$ with respect to $S$. One can show that the ring $R^{(-1)}R$ is absolutely flat. Clearly, $\eta(s)=\eta(s)^{2}(x_{s}+I)$ and $x_{s}+I=\eta(s)(x_{s}+I)^{2}$ where $\eta:R\rightarrow S^{(-1)}R$ is the canonical map and the pair $(S^{(-1)}R, \eta)$ satisfies in the following universal property: ``for each such pair $(A,\phi)$, i.e., $\phi:R\rightarrow A$ is a ring map and for each $s\in S$ there is some $c\in A$ such that $\phi(s)=\phi(s)^{2}c$ and $c=\phi(s)c^{2}$, then there exists a unique ring map $\psi:S^{(-1)}R\rightarrow A$ such that $\phi=\psi\circ\eta$''. It follows that if $\phi:A\rightarrow B$ is a ring map, then by the universal property there exists a unique ring map $\phi':=\phi^{(-1)}:A^{(-1)}A\rightarrow B^{(-1)}B$ such that the following diagram is commutative: $$\xymatrix{
A\ar[r]^{\phi}\ar[d]^{\eta_{1}} &B
\ar[d]^{\eta_{2}} \\A^{(-1)}A\ar[r]^{\phi'} &B^{(-1)}B}$$ where $\eta_{1}$ and $\eta_{2}$ are the canonical morphisms. The canonical ring map $\eta:R\rightarrow S^{(-1)}R$ induces a bijection between the corresponding prime spectra and $\Ker\eta\subseteq\mathfrak{N}$ where $\mathfrak{N}$ is the nil-radical of $R$. It is also easy to see that if $R$ is an absolutely flat ring, then the canonical ring map $\eta:R\rightarrow R^{(-1)}R$ is an isomorphism. For more details, see \cite{Olivier}. \\

Surjective ring maps are special cases of epimorphisms of rings. As a specific example, the canonical ring map $\mathbb{Z}\rightarrow\mathbb{Q}$ is an epimorphism of rings which is not surjective. Let $R$ be a ring. Then there exists a canonical injective flat epimorphism $\eta:R\rightarrow\mathcal{M}(R)$ such that every $R-$algebra whose structure morphism is an injective flat epimorphism can be canonically embedded in $\mathcal{M}(R)$. The ring $\mathcal{M}(R)$  is called ``the maximal flat epimorphic extension of $R$". For more details, see \cite[Proposition 3.4]{Lazard}.\\

We shall freely use the above facts in this article.\\

\section{Patch Closures}

In the following result, which we shall use later, the patch closures are computed in a certain way.\\

\begin{theorem}\label{lem 3} Let $R$ be a ring and let $E$ be a subset of $\Spec(R)$. Consider the patch topology on $\Spec(R)$ and also consider the canonical ring map $\pi: R\rightarrow\prod\limits_{\mathfrak{p}\in E}\kappa(\mathfrak{p})$, where $\kappa(\mathfrak{p})$ is the residue field of $R$ at $\mathfrak{p}$. Then $cl(E)=\Ima\pi^{\ast}$ where $cl(E)$ denotes the closure of $E$ in $\Spec(R)$ with respect to the patch topology.\\
\end{theorem}

{\bf Proof.} For each prime ideal $\mathfrak{p}\in E$, the canonical map $R\rightarrow\kappa(\mathfrak{p})$ factors as
$\xymatrix{R \ar[r]^{\pi} & A\ar[r]^{s_{\mathfrak{p}}} & \kappa(\mathfrak{p})}$ where $A:=\prod\limits_{\mathfrak{q}\in E}\kappa(\mathfrak{q})$ and $s_{\mathfrak{p}}$ is the canonical projection. Let $J_{\mathfrak{p}}:=s^{-1}_{\mathfrak{p}}(0)$. Then clearly $\mathfrak{p}=\pi^{-1}(J_{\mathfrak{q}})$. This implies that $cl(E)\subseteq\Ima\pi^{\ast}$. To prove the reverse inclusion we act as follows. First, we show that $cl(E')=\Spec(A)$ where $E'=\{J_{\mathfrak{p}} : \mathfrak{p}\in E\}$. Clearly the ring $A$ is absolutely flat. Therefore, by \cite[Theorem 3.3]{Abolfazl 2}, the Zariski and patch topologies over $\Spec(A)$ are the same. Let $f=(f_{\mathfrak{p}})_{\mathfrak{p}\in E}$ be a non-zero element of $A$. Thus there exists some $\mathfrak{p}\in E$ such that $f_{\mathfrak{p}}\neq0$. This implies that $J_{\mathfrak{p}}\in D(f)$.
Therefore, every non-empty standard Zariski open of $\Spec(A)$ meets $E'$. This means that $cl(E')=\Spec(A)$. We have $\pi^{\ast}(E')=E\subseteq cl(E)$. It follows that $E'\subseteq(\pi^{\ast})^{-1}\big(cl(E)\big)$. Thus, $cl(E')\subseteq(\pi^{\ast})^{-1}\big(cl(E)\big)$. Therefore, $\Ima\pi^{\ast}=\pi^{\ast}\big(\Spec(A)\big)=
\pi^{\ast}\big(cl(E')\big)\subseteq cl(E)$. $\Box$ \\

\section{Minimal and maximal spectra}

\begin{lemma} Let $\mathfrak{p}$ be a minimal prime of a ring $R$. If there is a Zariski quasi-compact open $U$ of $\Spec(R)$ such that $\mathfrak{p}\notin U$, then there exists an element $f\in R\setminus\mathfrak{p}$ such that $D(f)\cap U=\emptyset$. \\
\end{lemma}

{\bf Proof.} There are finitely many elements $g_{1},...,g_{n}\in R$ such that $U=\bigcup\limits_{i=1}^{n}D(g_{i})$. We have $g_{i}\in\mathfrak{p}$ for all $i$. It follows that the image of $g_{i}$ under the canonical map $R\rightarrow R_{\mathfrak{p}}$ is nilpotent. Thus, there is an element $f\in R\setminus\mathfrak{p}$ and a natural number $N$ such that $fg_{i}^{N}=0$ for all $i$. Then clearly $D(f)\cap U=\emptyset$. $\Box$ \\

\begin{proposition}\label{prop IVI} The minimal spectrum of a ring $R$ with respect to the Zariski topology is Hausdorff and totally disconnected. Moreover, $\Min(R)\cap D(f)$ is a clopen subset of $\Min(R)$ for all $f\in R$.\\
\end{proposition}

{\bf Proof.} It follows from the above Lemma. $\Box$ \\

\begin{lemma} Let $\mathfrak{m}$ be a maximal ideal of a ring $R$. If there is a flat quasi-compact open $U$ of $\Spec(R)$ such that $\mathfrak{m}\notin U$, then there exists an element  $f\in\mathfrak{m}$
such that $V(f)\cap U=\emptyset$.\\
\end{lemma}

{\bf Proof.} There is an element $g\in R$ such that $U\subseteq V(g)$ and $g\notin\mathfrak{m}$. Thus there are elements $a\in R$ and $f\in\mathfrak{m}$ such that $ag+f=1$. Clearly, $V(f)\cap V(g)=\emptyset$.  $\Box$ \\

As the dual of Proposition \ref{prop IVI} we have the following proposition. \\

\begin{proposition}\label{prop 1} The maximal spectrum of a ring $R$ with respect to the flat topology is Hausdorff and totally disconnected. Moreover, $\Max(R)\cap V(f)$ is a clopen subset of $\Max(R)$ for all $f\in R$.\\
\end{proposition}

{\bf Proof.} It follows from the above Lemma. $\Box$ \\

It is well known that the maximal spectrum of a ring is quasi-compact with respect to the Zariski topology. Dually, its minimal spectrum is quasi-compact with respect to the flat topology. But the minimal (respectively maximal) spectrum of a ring is not necessarily quasi-compact with respect to the Zariski (respectively flat) topology. As a specific example, $\Max(\mathbb{Z})$ is not quasi-compact with respect to the flat topology since the set of prime numbers is infinite. By \cite[Theorem 6 and Proposition 8]{Hochster 1}, there exists a ring $A$ whose prime ideals have precisely the reverse order of the primes of $\mathbb{Z}$ and so
$\Min(A)$ is homeomorphic to $\Max(\mathbb{Z})$. Therefore $\Min(A)$ is not Zariski quasi-compact. \\

\begin{theorem}\label{th 1} The maximal spectrum of a ring $R$ is compact with respect to the flat topology if and only if $R/J(R)$ is absolutely flat where $J(R)$ is the radical Jacobson of $R$. \\
\end{theorem}

{\bf Proof.} If $A:=R/J(R)$ is absolutely flat, then $\Max(A)=\Spec(A)$. Therefore, the image of the map $\Spec(A)\rightarrow\Spec(R)$ induced by the canonical map $R\rightarrow A$ is equal to $\Max(R)$. Hence, $\Max(R)$ is quasi-compact. Conversely, let $f\in R\setminus J(R)$. Let $\mathfrak{m}$ be a maximal ideal of $R$ such that $f\notin\mathfrak{m}$. Then there are elements $a\in R$ and $c\in\mathfrak{m}$ such that $1=af+c$. Now by applying the quasi-compactness of $\Max(R)$ then we may find a finite number of elements $a_{1},...,a_{n}$ and $c_{1},...,c_{n}$ of $R$ such that $\Max(R)\subseteq\V(f)\cup\big(\bigcup\limits_{i=1}^{n}\V(c_{i})\big)$ and $1=a_{i}f+c_{i}$ for all $i$. This implies that $fc_{1}...c_{n}\in J(R)$ and $1=bf+c_{1}...c_{n}$ for some $b\in R$. It follows that $f-bf^{2}\in J(R)$. Therefore, $R/J(R)$ is absolutely flat. $\Box$ \\

Theorem \ref{th 1} implies that $\Max(R)$ is compact in the flat topology iff $R/J(R)$ has Krull dimension zero. \\

In order to establish the dual of Theorem \ref{th 1}, the following results are required.\\

\begin{lemma}\label{lemm 5} Let $R\subseteq S$ be an extension of rings with $R$ absolutely flat then $S$ is $R-$faithfully flat.\\
\end{lemma}

{\bf Proof.} If $\mathfrak{m}$ is a maximal ideal of $R$, then it is a minimal prime of $R$. So there exists a (minimal) prime $\mathfrak{p}$ of $S$ lying over it. It follows that $\mathfrak{m}S\subseteq\mathfrak{p}\neq S$. Hence, $S$ is $R-$faithfully flat. $\Box$ \\

\begin{lemma}\label{lemma 9} Every injective flat ring map $\phi:R\rightarrow A$ with $A$ absolutely flat factors through an injective flat epimorphism $\psi:R\rightarrow B$ followed be the canonical injection $i:B\rightarrow A$ such that $B$ is absolutely flat.\\
\end{lemma}

{\bf Proof.} Consider the following commutative diagram:
$$\xymatrix{
R\ar[r]^{\phi}\ar[d]^{\eta_{1}}&A \ar[d]^{\eta_{2}}\\R^{(-1)}R\ar[r]^{\phi'} & A^{(-1)}A}$$
where $R^{(-1)}R$ and $A^{(-1)}A$ are the pointwise localizations of $R$ and $A$, respectively. By the hypotheses, $\eta_{2}$ is an isomorphism. Let $\theta:=\eta_{2}^{-1}\circ\phi':R^{(-1)}R\rightarrow A$ also let $B:=\Ima\theta$ which is an absolutely flat ring since $R^{(-1)}R$ is absolutely flat. The map $\theta$ factors as $i\circ\theta'$ where $\theta':R^{(-1)}R\rightarrow B$ is the canonical surjection which is induced by $\theta$. Let $\psi:=\theta'\circ\eta_{1}:R\rightarrow B$ which is an epimorphism.
By lemma \ref{lemm 5}, $i$ is faithfully flat. Therefore, $\psi$ is injective and flat since $\phi=i\circ\psi$. $\Box$ \\

\begin{lemma}\label{lemma 7} Let $\phi:R\rightarrow A$ be a ring map such that for each prime $\mathfrak{q}$ of $A$,  $A_{\mathfrak{q}}$ is $R_{\mathfrak{p}}$-flat where $\mathfrak{p}:=\phi^{-1}(\mathfrak{q})$. Then $A$ is $R-$flat.\\
\end{lemma}

{\bf Proof.} It is an easy exercise. $\Box$ \\

\begin{theorem}\label{th 2}\cite[Proposition 1]{Quentel} If the minimal spectrum of a reduced ring $R$ is Zariski compact, then the ring $\mathcal{M}(R)$, the maximal flat epimorphic extension of $R$, is absolutely flat.\\
\end{theorem}

{\bf Proof.} The minimal spectrum of $R$ is a patch closed subset of $\Spec(R)$. Because suppose there is some $\mathfrak{q}\in cl\big(\Min(R)\big)$ which is not in $\Min(R)$.
Then for each $\mathfrak{p}\in\Min(R)$ there exists some $f\in\mathfrak{q}$ which is not in $\mathfrak{p}$. By applying the quasi-compactness of $\Min(R)$ then we may find finitely many elements $f_{1},...,f_{n}\in\mathfrak{q}$ such that $\Min(R)\subseteq\bigcup\limits_{i=1}^{n}D(f_{i})$. Clearly $\bigcap\limits_{i=1}^{n}V(f_{i})$ is a patch open neighbourhood of $\mathfrak{q}$. Therefore it meets $\Min(R)$. But this is a contradiction.
Thus, $cl\big(\Min(R)\big)=\Min(R)$. Let $A:=\prod\limits_{\mathfrak{p}\in\Min(R)}R_{\mathfrak{p}}$. It is absolutely flat since $R_{\mathfrak{p}}=\kappa(\mathfrak{p})$ for all $\mathfrak{p}\in\Min(R)$. By Theorem \ref{lem 3}, $\Min(R)=\Ima\pi^{\ast}$ where $\pi:R\rightarrow A$ is the canonical map.  Hence, by Lemma \ref{lemma 7}, $\pi$ is flat. It is also injective. Thus, by Lemma \ref{lemma 9}, there exists an injective flat epimorphism $\psi:R\rightarrow B$ such that $B$ is absolutely flat. By \cite[Proposition 3.4]{Lazard}, there is a (unique injective) ring map $h:B\rightarrow S$ such that $\eta=h\circ\psi$ where
$S:=\mathcal{M}(R)$ and $\eta:R\rightarrow S$ is the canonical map. It follows that for each prime $\mathfrak{q}$ of $S$ then $\mathfrak{p}:=\eta^{-1}(\mathfrak{q})$ is a minimal prime of $R$ because $B$ is absolutely flat and every flat ring map has the going-down property. The canonical map $\kappa(\mathfrak{p})=R_{\mathfrak{p}}\rightarrow R_{\mathfrak{p}}\otimes_{R}S_{\mathfrak{q}}
\simeq S_{\mathfrak{q}}$ is a flat epimorphism since flat maps and epimorphisms are stable under base change. It is also faithfully flat since $\mathfrak{p}S\neq S$. Every faithfully flat epimorphism is an isomorphism. Therefore, $S$ is absolutely flat. $\Box$ \\

\begin{lemma}\label{lemma 1} Let $\phi:R\rightarrow A$ be an injective flat ring map. If $\Min(A)$ is Zariski compact then $\Min(R)$ is as well. \\
\end{lemma}

{\bf Proof.} It $\mathfrak{p}$ is a minimal prime of $R$, then there exists a minimal prime of $A$ lying over $\mathfrak{p}$ since $\phi$ is injective. Therefore $\Min(R)\subseteq\phi^{\ast}\big(\Min(A)\big)$. Moreover, the inverse image of every minimal prime of $A$ under $\phi$ is a minimal prime of $R$ since every flat ring map has the going-down property. Therefore, $\phi^{\ast}\big(\Min(A)\big)=\Min(R)$. Hence $\Min(R)$ is quasi-compact. $\Box$ \\

Now as the dual of Theorem \ref{th 1} we have the following corollary. \\

\begin{corollary}\label{coro 2}\cite[Proposition 1]{Quentel} The minimal spectrum of a ring $R$ is Zariski compact if and only if $\mathcal{M}(R/\mathfrak{N})$ is absolutely flat where $\mathfrak{N}$ is the nil-radical of $R$.\\
\end{corollary}

{\bf Proof.} The canonical map $R\rightarrow R/\mathfrak{N}$ induces a homeomorphism between the corresponding prime spectra which maps $\Min(R/\mathfrak{N})$ onto $\Min(R)$. Then apply Theorem \ref{th 2} for the implication ``$\Rightarrow$" and Lemma \ref{lemma 1} for the reverse. $\Box$ \\

\begin{corollary}\label{Corollary I} The minimal spectrum of a reduced ring $R$ is Zariski compact if and only if the total ring of fractions of the polynomial ring $R[x]$ is absolutely flat.\\
\end{corollary}

{\bf Proof.} The ring $R[x]$ is reduced and every minimal prime of it is of the form $\mathfrak{p}[x]$ where $\mathfrak{p}$ is a minimal prime of $R$. Suppose $\Min(R[x])\subseteq\bigcup\limits_{i}D(f_{i})$ where $f_{i}\in R[x]$ for all $i$. Then $\Min(R)\subseteq\bigcup\limits_{i}D(c_{i})$ where $c_{i}$ is some coefficient of $f_{i}$. Thus, the compactness of the minimal spectrum of $R$ is equivalent to that of $R[x]$. Then apply Corollary \ref{coro 2} and the fact that the total ring of fractions of the polynomial ring $R[x]$ is canonically isomorphic to $\mathcal{M}(R[x])$. $\Box$ \\

\begin{remark} Note that ``reduced ring'' assumption in Corollary \ref{Corollary I} is crucial, because without this assumption then the assertion does not necessarily hold. As an example, take $R=\mathbb{Z}/4\mathbb{Z}$ then $\Min(R)$ is compact since it is finite, but $T(R[x])$ is not absolutely flat, since every absolutely flat ring and so each subring are reduced but $R$ is not reduced. \\
\end{remark}

\section{Noetherian property}

Here, we give a characterization for noetherianess of the Zariski topology in terms of the flat topology.\\

\begin{theorem}\label{th 3} The prime spectrum of a ring $R$ is noetherian with respect to the Zariski topology if and only if the flat opens of $\Spec(R)$ are stable under the arbitrary intersections.\\
\end{theorem}

{\bf Proof.} First assume that $\Spec(R)$ is Zariski noetherian. If $I$ is an ideal of $R$ then $V(I)$ is a flat open, because $U=\Spec(R)\setminus\V(I)$ is Zariski quasi-compact and hence there exists a finite set $\{f_{1},...,f_{n}\}$ of elements of $R$ such that  $U=\bigcup\limits_{i=1}^{n}D(f_{i})$. It follows that $\V(I)=\bigcap\limits_{i=1}^{n}\V(f_{i})$. To prove the assertion it suffices to show that the intersection of every family of basis flat opens is flat open.
The basis flat opens are precisely of the form $V(I)$ where $I$ is a finitely generated ideal of $R$. Let $\{I_{\alpha}\}$ be a family of ideals of $R$. If $\mathfrak{p}\in\bigcap\limits_{\alpha}\V(I_{\alpha})$ then clearly $\V(\mathfrak{p})\subseteq\bigcap\limits_{\alpha}\V(I_{\alpha})$. Therefore, $\bigcap\limits_{\alpha}\V(I_{\alpha})$ is flat open. To prove the reverse, it suffices to show that every Zariski open $U=\Spec(R)\setminus\V(I)$ of $\Spec(R)$ is quasi-compact where $I$ is an ideal of $R$. But $V(I)$ is flat open since $\V(I)=\bigcap\limits_{f\in I}\V(f)$. It is also quasi-compact. Thus, there exists a finite set $\{I_{1},...,I_{s}\}$ of finitely generated ideals of $R$ such that $\V(I)=\bigcup\limits_{k=1}^{s}\V(I_{k})=\V(I_{1}I_{2}...I_{s})$. Note that $I_{1}I_{2}...I_{s}=(f_{1},...,f_{n})$ is a finitely generated ideal of $R$. Therefore, $\sqrt{I}$ is equal to the radical of a finitely generated ideal of $R$. We have then $U=\Spec(R)\setminus\V(\sqrt{I})=\Spec(R)\setminus\bigcap\limits_{i=1}^{n}
V(f_{i})=\bigcup\limits_{i=1}^{n}D(f_{i})$.  $\Box$ \\

\textbf{Acknowledgements.} The author would like to give sincere thanks to Professor Irena Swanson for valuable comments about the paper. The author would also like to present heartfelt thanks to the referee for careful reading of the paper and for many valuable comments which improved the paper. \\

\end{document}